\documentclass[a4paper, 10pt, conference]{ieeeconf}      

\IEEEoverridecommandlockouts                              
\usepackage{graphics} 
\usepackage{epsfig} 
\usepackage{mathptmx} 
\usepackage{times} 
\usepackage{amsmath} 
\usepackage{amssymb}  
\graphicspath{{images/}}

\title{\LARGE \bf
Centralized Recursive Optimal Scheduling\\ of Parallel Buck Regulated Battery Modules
}

\author{Yunfeng Jiang, Abdulelah H. Habib, Xin Zhao, Louis J. Shrinkle, Raymond A. de Callafon
        \thanks{Yunfeng Jiang, Abdulelah H. Habib and Xin Zhao are with the Department of Mechanical and Aerospace Engineering,  University of California San Diego, 9500 Gilman Drive, Mail Code 0411, La Jolla, CA 92093-0411, USA
        {\tt\small yuj034@eng.ucsd.edu, ahhabib@eng.ucsd.edu, xiz028@ucsd.edu}}%
\thanks{Raymond A. de Callafon is with Faculty of Mechanical and Aerospace Engineering, University of California San Diego, 9500 Gilman Drive, Mail Code 0411, La Jolla, CA 92093-0411, USA
        {\tt\small callafon@ucsd.edu}}%
        \thanks{Louis J. Shrinkle is with 1856 Wilstone Ave, Encinitas, CA 92024, USA}
}

\begin{document}

\maketitle
\thispagestyle{empty}
\pagestyle{empty}

\begin{abstract}

This paper presents a centralized recursive optimal scheduling method for a battery system that consists of parallel connected battery modules with different open circuit voltages and battery impedance characteristics. Examples of such a battery system can be found in second-life, exchangeable or repurposed battery systems in which batteries with different charge or age characteristics are combined to create a larger storage capacity. The proposed method in this paper takes advantage of the availability of buck regulators in the battery management system (BMS) to compute the optimal voltage adjustment of the individual modules to minimize the effect of stray currents between the parallel connected battery modules. Our proposed method recursively computes the optimal current scheduling that balances (equals) each module current and maximize total bus current without violating any of the battery modules operating constraints. Recursive implementation guarantees robust operation as the battery modules operating parameters change as the battery pack (dis)charges or ages. In order to demonstrate the capability of this method in real battery system, an experimental setup of 3 parallel placed battery modules is built. The experimental results validate the feasibility and show the advantages of this current scheduling method in a real battery application, despite the fact that each module may have different impedance, open circuit voltage and charge parameters. 
\end{abstract}

\section{INTRODUCTION}
Electric vehicles (EVs) are widely regarded as a promising environmental-friendly solution for future automotive industry, due to technological developments and an increased focus on renewable energy \cite{xing_battery_2011,chaturvedi_algorithms_2010,tribioli_energy_2016}. Most EVs use lithium-ion batteries (LIBs) to storage energy and supply power to the electric grid including communication and control systems, because LIBs have higher specific power and energy density, longer life span and lower self-discharge rate than most other practical batteries \cite{li_method_2016,cai_impact_2017}. 

Typically, a series connection of LIBs are used to create a battery module that achieves a desired open circuit or battery terminal voltage, while a parallel connection of battery modules is used to increase the total energy storage capacity of the battery pack. Parallel connected battery modules also increase power capacity to both fulfill an acceptable driving range and maintain high performance of EVs \cite{khaligh_battery_2010,liu_sensor_2017,finesso_cost-optimized_2016,chen_model-based_2016}. Such configuration also shows significant benefits in second-life battery applications, such as demand charge management, renewable energy integration and regulation energy management in EVs \cite{deng_life_2017,tong_off-grid_2013}.

Unfortunately, variability in the production process of LIBs can not ensure identical battery modules and moreover, parallel placed battery modules exposed to the same (dis)charge cycles might age differently. Both effects cause discrepancies in internal resistance (impedance), temperature gradients, and operation conditions, such as power storage/delivery demand and environmental conditions. These discrepancies may limits the ability to extract or store the full electrical energy capacity in the battery system \cite{lu_review_2013,wieczorek_mathematical_2017,wang_novel_2015,sarwar_experimental_2016}. Therefore, it is essential to build a battery management system (BMS) for a high power battery pack with parallel placed battery modules, so as to accurately estimate the state of charge (SOC) and state of health (SOH) of the modules to protect the battery from operating outside its Safe Operating Area (SOA) \cite{jiang_data-based_2017, jiang_dynamic_2016,zhao_modeling_2016,zhao_data-based_2013-1,xia_accurate_2016}. 

\begin{figure}[ht]
\begin{center}
\includegraphics[width=\columnwidth]{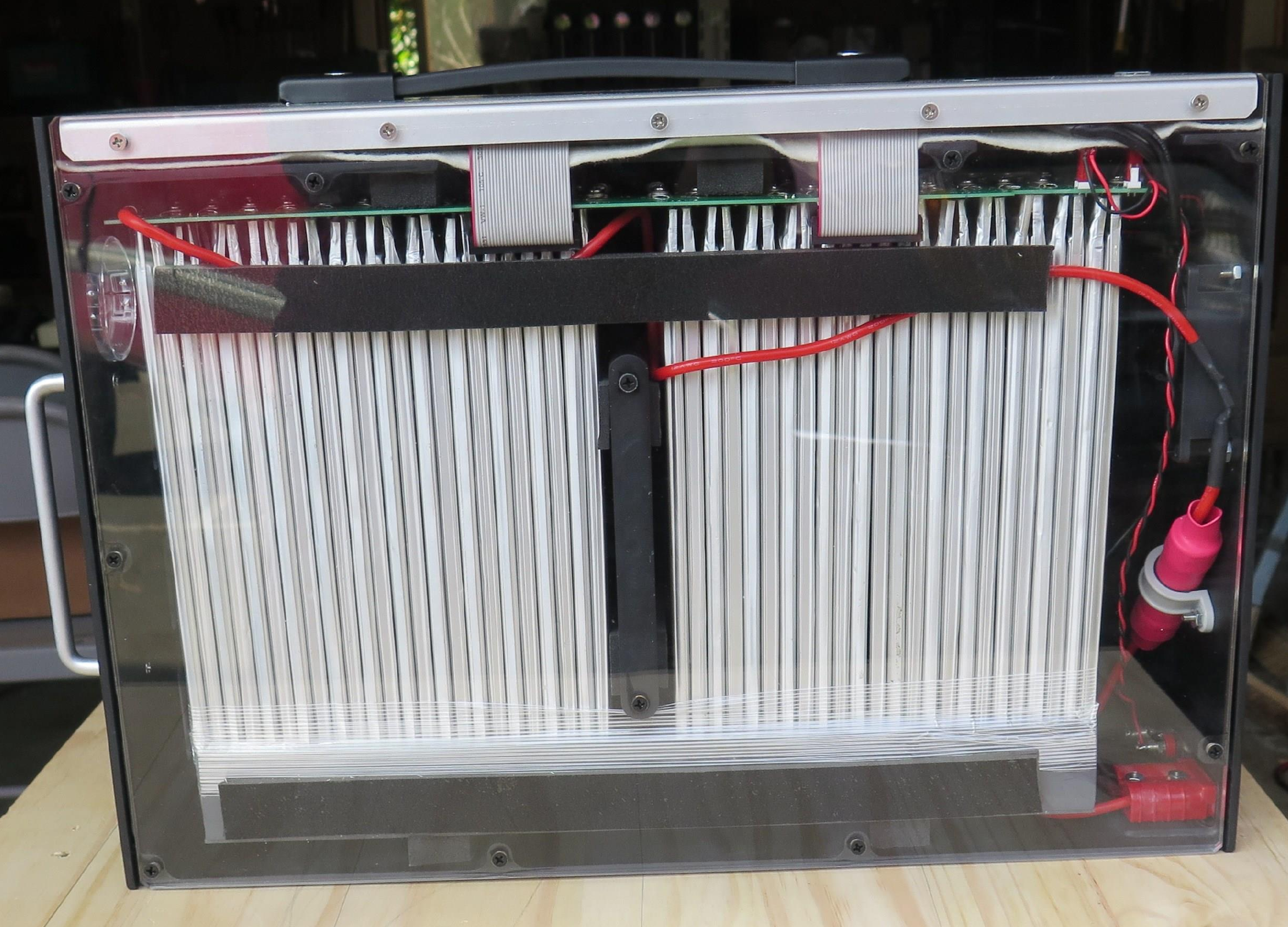}
\caption{Exchangeable battery module with a series connection of LIBs in a suitcase size format. Multiple of these battery modules are connected in parallel to increase power and energy storage capabilities.} 
\label{fig:exchangeable battery module with parallel connected battery cells}
\end{center}
\end{figure}

A high power battery pack with parallel connected battery modules that allow exchangeable modules can be viable alternative to increase power storage capacity and operation efficiency in a flexible way \cite{zhao_current_2014,brand_current_2016}. In this paper we consider the use of separate and exchangeable battery modules that consists of a series connection of LIBs in a suitcase format, as shown in Fig.~\ref{fig:exchangeable battery module with parallel connected battery cells}. Several LIBs are connected in series to achieve a large open circuit voltage, whereas a buck regulator is included in the battery module to control the battery module terminal voltage. The buck regulator uses a digital Pulse Width Modulation (PWM) signal to several MOSFETs with an RLC filter (an electrical circuit consisting of a resistor (R), an inductor (L), and a capacitor (C)) to effectively reduce the full scale OCV. Subsequently, multiple modules can be placed in parallel on a single DC connection bus to create the full battery pack/system with a desired power and energy storage specification.

With such configuration, the efficiency and flexibility of this particular battery system is significantly improved, because partly empty or failing modules of the battery pack can be exchanged for fast charging and fault correction capabilities. As such, this configuration is applicable to second-life battery applications and/or EVs with partly exchangeable batteries. However, there are some challenges and bottlenecks with this configuration: each individual battery module may have different SOCs, instantaneous and nominal capacity, and internal impedances. Therefore, the terminal voltage of each parallel placed battery module must be controlled, and denoted by 'scheduling' in this paper, when charging and discharging battery modules.

The scheduling problem for parallel placed battery modules is solved in this paper by computing the optimal voltage adjustment of the individual modules to minimize the effect of stray currents between the parallel connected battery modules. The scheduling can compute various charging/discharging solutions and we consider a solution where the module currents are closely matched. The actual implementation of voltage regulation is achieved via the digital PWM signal used in the buck regulator included in the battery module. The scheduling solution proposed in this paper computes the PWM in \% of full scale for each module, without the explicit knowledge of the electrical module parameters of open circuit voltage and impedance of each module and the impedance/load connected to the battery pack. A recursive formulation of the scheduling solution allows estimation of electrical module parameters and adaptation to time varying load conditions on the battery pack. Finally, experimental results are included in this paper to validate the feasibility and performance of the proposed current scheduling control technique in different scenarios of a real battery system. 

\section{PARALLEL BUCK REGULATED BATTERY MODULES}

\subsection{Introduction and Assumptions}

For scheduling of parallel placed battery modules, we assume that each module $k$ is characterized by an modulated ideal voltage supply in series with an impedance. For each module $k$ we make the following assumptions:

\begin{itemize}

\item The ideal voltage supply is given by $V_k = \alpha_k V^{OCV}_k$, where $V^{OCV}_k$ is the open circuit voltage (OCV) or terminal voltage of the battery module in case of no load and where voltage modulation can be only down via $0 < \alpha_k < 1$. 

\item The internal impedance of a module is given by $Z_k$ and is slowly time varying.
\end{itemize}

\begin{figure}
\begin{center}
\includegraphics[width=\columnwidth]{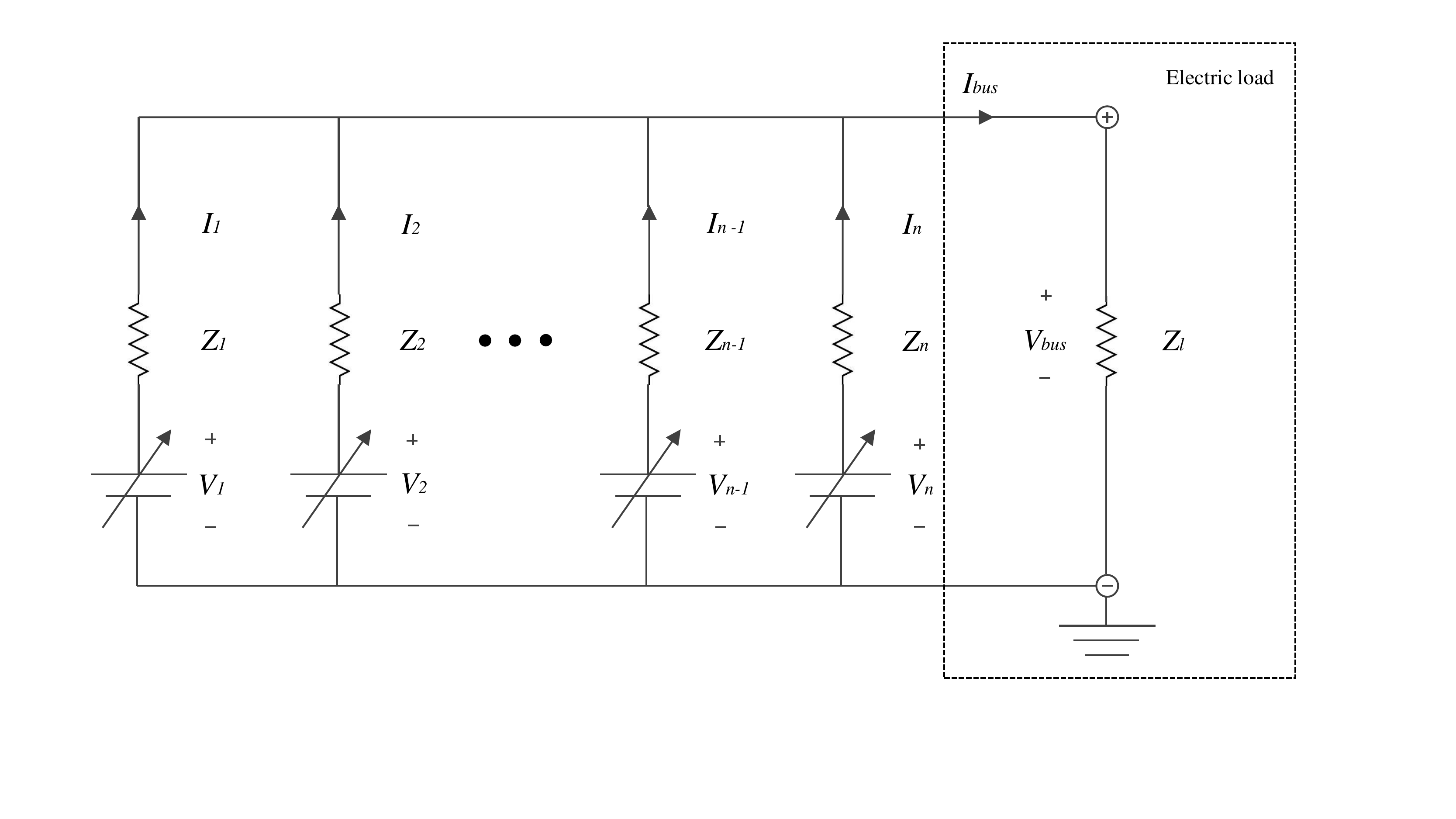}
\caption{Model for current scheduling.} 
\label{fig:model for current scheduling}
\end{center}
\end{figure}

The slowly time varying nature of the model impedance $Z_k$ is in comparison to the time varying nature of the external load impedance $Z_l$ connected to the battery pack, as shown in Fig.~\ref{fig:model for current scheduling}. Following the electrical diagram of Fig.~\ref{fig:model for current scheduling}, application of Kirchhoff's current and voltage law now leads to the following results for parallel placed modules characterized by an modulated ideal voltage supply $V_k$ in series with an impedance $Z_k$:

\begin{itemize}

\item The current $I_k$ of each module ensures that the bus current
\begin{equation}
I_{bus} = \sum_{k=1}^n I_k
\label{eq:Ibus}
\end{equation}
due to Kirchhoff's current law

\item The bus voltage $V_{bus}$ satisfies
\begin{equation}
V_{bus} = V_k - Z_k I_k 
\label{eq:Vbus}
\end{equation}
for each module $k$ due to Kirchhoff's voltage law.
\end{itemize}

\subsection{Bus voltage as a function of module voltages}

The above results can be combined to compute the bus voltage $V_{bus}$ or bus current $I_{bus}$, when a load $Z_l$ is applied to battery pack that consists of a set of parallel placed modules. For a given set of values for the modulated voltages $V_k,~k=1,2,\ldots,n$, with the individual module currents $I_k$ given by
\begin{equation}
I_k = \frac{V_k- V_{bus}}{Z_k}
\label{eq:Ik}
\end{equation}
we have
\[
V_{bus}= Z_{l} I_{bus} = Z_{l} \sum_{k=1}^n I_k = Z_{l} \sum_{k=1}^n \frac{V_k- V_{bus}}{Z_k}
\]
From this last expression we can solve $V_{bus}$ via
\[
V_{bus} = Z_{l} \sum_{k=1}^n \frac{V_k}{Z_k} - Z_{l} V_{bus} \sum_{k=1}^n \frac{1}{Z_k}
\]
making
\[
 V_{bus}= \frac{\displaystyle\sum_{k=1}^n \frac{V_k}{Z_k}}{\displaystyle\frac{1}{Z_l} + \sum_{k=1}^n \frac{1}{Z_k}}
\]
Although this expressions may look complicated at first, it is important to recognize that the bus voltage $V_{bus}$ is determined by the linear combination
\begin{equation}
\begin{array}{rcl}
V_{bus} &=& g_1 V_1 + g_2 V_2 + \ldots g_n V_n~~~ \mbox{where}\\
g_j &=& \displaystyle \frac{\displaystyle\frac{1}{Z_j}}{\displaystyle\frac{1}{Z_l} + \sum_{k=1}^n \frac{1}{Z_k}},~~ j=1,2,\ldots,n
\end{array}
\label{eq:Vgains}
\end{equation}
in which the "gain factors" $g_k,~ k=1,2,\ldots,n$ are given by a familiar parallel connection formula of impedances $Z_k$ and $Z_l$. 

\subsection{Module currents as a function of module voltages}

With knowledge of $V_{bus}$, clearly also $I_{bus} = V_{bus}/Z_l$ can be computed. It is tempting to write $I_{bus}$ also as a similar linear combination by using (\ref{eq:Vgains}) to obtain
\[
I_{bus} = \frac{V_{bus}}{Z_l} = \frac{g_1}{Z_l} V_1 + \frac{g_2}{Z_l} V_2 + \ldots + \frac{g_n}{Z_l} V_n
\]
and conclude based on (\ref{eq:Ibus}) that the individual module currents are given by $I_k = \frac{g_k}{Z_l} V_k$, but that is incorrect. The individual module currents $I_k$ are typically a linear combination of {\em all\/} the modulated module voltages $V_k$, for which an expression is derived here.

With the individual module currents $I_k$ given in (\ref{eq:Ik}) and $V_{bus}$ given in (\ref{eq:Vbus}), we can then obtain
\[
I_k = \frac{1}{Z_k} \left ( \frac{\displaystyle\frac{V_k}{Z_l} + \sum_{m=1}^n \frac{V_k}{Z_m} - \displaystyle\sum_{m=1}^n \frac{V_m}{Z_m}}{\displaystyle\frac{1}{Z_l} + \sum_{m=1}^n \frac{1}{Z_m}} \right )
\]
where the summation index has been changed to $m$ to avoid confusion with the specific module current $I_k$ indexed with $k$.
The expression for $I_k$ can be simplified to the insightful expression
\begin{equation}
\begin{array}{rcl}
I_k &=& \displaystyle d_{k,1} V_1 + d_{k,2} V_2 + \ldots + d_{k,n} V_n,~~~ \mbox{where}\\ 
d_{k,j} &=& \left \{ \begin{array}{rl}\displaystyle -\frac{1}{Z_k} \cdot \frac{\displaystyle\frac{1}{Z_j}}{\displaystyle\frac{1}{Z_l} + \sum_{m=1}^n \frac{1}{Z_m}} & \mbox{for} ~ j \not = k\\
\displaystyle \frac{1}{Z_k} \cdot \frac{\displaystyle\frac{1}{Z_l} + \sum_{m=1}^{n} \frac{1}{Z_m} - 
\frac{1}{Z_k}}{\displaystyle\frac{1}{Z_l} + \sum_{m=1}^n \frac{1}{Z_m}} & \mbox{for} ~ j = k
\end{array} \right .
\end{array}
\label{eq:Igains}
\end{equation}
where the coefficients $d_{k,j},~k=1,2,\ldots,n$ and $j=1,2,\ldots,n$ build up a $n \times n$ impedance matrix $D=[d_{k,j}]$. The impedance matrix $D=[d_{k,j}]$ relates module currents $I_k$ to module voltage $V_k$ according to
\begin{equation}
\left [ \begin{array}{c} I_1 \\ I_2 \\ \vdots \\ I_n \end{array} \right ] = \left [ \begin{array}{cccc} d_{1,1} & d_{1,2} & \cdots & d_{1,n} \\ d_{2,1} & d_{2,2} & \cdots & d_{2,n} \\ 
\vdots & \vdots & \cdots & \vdots \\
d_{n,1} & d_{n,2} & \cdots & d_{n,n} \end{array} \right ]
\left [ \begin{array}{c} V_1 \\ V_2 \\ \vdots \\ V_n \end{array} \right ]
\label{eq:impedancematrix}
\end{equation}
with $d_{k,j}$ given in (\ref{eq:Igains}), and will be useful for the explicit computation of module currents $I_k$ as a function of the module voltages $V_k$ and {\em visa versa}. 

It can be observed from the definition of the matrix $D=[d_{k,j}]$ that $D$ is symmetric. With all resistive values positive, it can also be shown that $D$ is also positive definite, making $D$ non-singular. With $D$ invertible, we can also compute module voltages $V_k$ as a function of desired module currents $I_k$ for the parallel placed battery modules.

\section{CENTRALIZED RECURSIVE OPTIMAL CURRENT SCHEDULING}

\subsection{Relative scaling of module currents}

Given the knowledge on the internal impedances $Z_k,~ k=1,2,\ldots,n$ and a fixed (but unknown) load impedance $Z_l$, the idea of module scheduling is formulated in this paper as the following problem. Compute the buck regulated module voltage $V_k \leq V_k^{OCV}$, such that module currents $I=[I_1~I_2~ \cdots I_n]^T$ are scaled to satisfy
\begin{equation}
I = \beta \left [ \begin{array}{c} \beta_1 \\ \beta_2 \\ \vdots \\ \beta_n \end{array} \right ],~~ 0 \leq \beta_k \leq 1,~ k=1,2,\ldots,n
\label{eq:Ibal}
\end{equation}
in which $\beta$ is used for absolute scaling, whereas $0 \leq \beta_k \leq 1$ specifies the relative scaling of the module current $I_k$. The value $\beta$ satisfies $\beta>0$ for battery module discharging, whereas $\beta<0$ for battery charging. Ideally, module scheduling should be done {\em despite\/} the lack of knowledge on the internal module impedance $Z_k$ and the externally applied load impedance $Z_l$. A recursive solution will be formulated to accomplish this later in the paper.

The motivation for the relative scaling $\beta_k$ of module currents $I$ according to (\ref{eq:Ibal}) is to discharge/charge current $I_k$ out/in of a module $k$ based on the individual SOC of each battery module. For example, if the SOC of module $k$ is denoted by $0\%< SOC_k\leq100\%$, an appropriate relative scaling $\beta_k$ of a module current $k$ could be defined as
\[
\beta_k = \frac{SOC_k}{\max_{k=1,2,\ldots,n} SOC_k} \leq 1
\]
to ensure that battery modules with a smaller SOC will discharge less current compared to battery modules with a larger SOC. Similarly, for charging we may want
\[
\beta_k = \frac{\min_{k=1,2,\ldots,n} SOC_k}{SOC_k} \leq 1
\]
to satisfy that battery modules with a smaller SOC will charge faster with a larger current compared to battery modules with a larger SOC. If all modules have the same storage capacity with the same relative SOC and are required to follow the same (dis)charging profile, the relative scaling $\beta_k$ of module currents $I$ according to (\ref{eq:Ibal}) can be required to satisfy $\beta_k = 1,~ k=1,\ldots,n$ causing 
\begin{equation}
I_1 = I_{2} = \cdots = I_n
\label{eq:equal}
\end{equation}
will be denoted by {\em equal SOC balancing\/} in this paper.

\subsection{Module scheduling via Linear Programming}

In case of full information on the external load impedance $Z_l$ and the internal impedance $Z_k$, the optimal modulated module voltages $V_k$ can be computed directly.
With the definition of the (invertible) impedance matrix $D$ in (\ref{eq:impedancematrix}) we can actually compute the set of internal module voltages $V=[V_1~V_2~\cdots~V_n]^T$ directly from a desired set of module currents $I=[I_1~I_2~\cdots~I_n]^T$. Using the vector notation
\[
V = \left [ \begin{array}{c} V_1 \\ V_2 \\ \vdots \\ V_n \end{array} \right ],~~ V^{OCV} = \left [ \begin{array}{c} V_1^{OCV} \\ V_2^{OCV} \\ \vdots \\ V_n^{OCV} \end{array} \right ]
\]
and the module currents in the vector format $I$ of (\ref{eq:Ibal}), the problem of module scheduling requires the computation of the maximum value of the current scaling $\beta>0$ such that $V \leq V^{OCV}$.

With the (invertible) impedance matrix $D$ in (\ref{eq:impedancematrix}), the problem of module scheduling can be written as a linear programming (LP) problem that can compute a globally optimal value of the current values $\beta$. By recognizing that
\[
V = D^{-1} \left [ \begin{array}{cccc} \beta_1 & \beta_2 & \cdots & \beta_n \end{array} \right ]^T \beta
\]
and the optimization for $\beta\geq 0$ (for discharging) of the module currents can be written as
\[
\begin{array}{c}
\displaystyle 
\max_\beta \beta \\
\mbox{s.t.}~  D^{-1} \left [ \begin{array}{cccc} \beta_1 & \beta_2 & \cdots & \beta_n \end{array} \right ]^T \beta \leq V^{OCV}
\end{array}
\]
and equivalent to a LP problem
\begin{equation}
\begin{array}{c}
\displaystyle \beta_{opt} = \min_\beta f^T \beta,~ \mbox{s.t.}~ A\beta \leq b,~~ \mbox{with}\\ ~~ A=D^{-1} \left [ \begin{array}{cccc} \beta_1 & \beta_2 & \cdots & \beta_n \end{array} \right ]^T,~~ \\ f^T = - 1,~~ \mbox{and}~~ b=V^{OCV}
\end{array}
\label{eq:linprog}
\end{equation}
The LP problem in (\ref{eq:linprog}) will compute the optimal value of $\beta_{opt}$. Once $\beta_{opt}$ is know, the (optimal) module
currents are given by
\[
I_{opt} = \beta_{opt} \left [ \begin{array}{cccc} \beta_1 & \beta_2 & \cdots & \beta_n \end{array} \right ]^T
\]
and the (optimal) module voltages are readily computed via
\begin{equation}
V_{opt} = \beta_{opt} D \left[ \begin{array}{cccc} \beta_1 & \beta_2 & \cdots & \beta_n \end{array} \right ]^T
\label{eq:Vopt}
\end{equation}
based on (\ref{eq:impedancematrix}).

\subsection{Centralized recursive module scheduling}

The LP solution in (\ref{eq:linprog}) requires knowledge of the impedance matrix $D$ in (\ref{eq:impedancematrix}) that is fully characterized by the internal module impedances $Z_k$ and the external load impedance $Z_l$. Once the impedance matrix $D$ is known, the linear optimization problem in (\ref{eq:linprog}) can be solved to compute the optimal value $\beta_{opt}$ of the (equal) balancing module currents, given the constraints on the OCV's $V^{OCV}$ for each module.

Although the internal module impedances $Z_k$ may be monitored by the battery management system (BMS), the external load $Z_l$ may not be known. In fact, the external load $Z_l$ may be time-varying fairly fast due to varying power demands, while internal module impedances $Z_k,~ k=1,2,\ldots,n$ typically only vary slowly over time. Clearly, module scheduling must be done without explicit knowledge of the external load $Z_l$. In recursive module scheduling we will update the impedance matrix $D$ recursively to allow for the computation of the optimal modulated module voltages $V_k$. The external load $Z_l$ can be estimated by monitoring the bus voltage $V_{bus}$ {\em and\/} the bus current $I_{bus}$. Using Ohm's law, we may estimate $Z_l = \frac{V_{bus}}{I_{bus}}$ so that the value of $Z_l$ in the impedance matrix $D$ can be replaced by the ratio of $V_{bus}$ and $I_{bus}$. 

Since the optimal values of the internal module voltages $V$ and the resulting bus voltage $V_{bus}$ again depend on the impedance matrix $D$, a straightforward recursive procedure can be used to recursively update $D$ and compute the optimal modulated module voltages $V$ for (scaled) balancing module currents. Starting from an initial choice for the internal module voltages $V=[V_1~V_2~\cdots~V_n]^T$,  the bus voltage $V_{bus}$ and the bus current $I_{bus}$ are measured to compute $Z_l = \frac{V_{bus}}{I_{bus}}$. With knowledge of the external load $Z_l$, the impedance matrix $D$ in (\ref{eq:impedancematrix}) can be updated and used to solve the LP problem in (\ref{eq:linprog}) to obtain the optimal current scaling $\beta_{opt}$. With $\beta_{opt}$ and the the impedance matrix $D$, the resulting (optimal) individual module voltages $V_{opt}$ in (\ref{eq:Vopt}) can be communicated to each of the modules. The procedure can be implemented recursively in time and summarized in the procedure below.

\textbf{\em Procedure:\/}  Assuming fixed internal impedances $Z_k,~ k=1,2,\ldots,n$ but a time-varying load impedance $Z_l$, the (centralized) recursive implementation is as follows:

\begin{enumerate}

\item Set time index $t=0$ and communicate the $n$ elements $V_k[0]$ of the initial module voltages $V[0] = [V_1[0]~V_2[0]~ \cdots~ V_n[0]]^T$ to each of the modules $k=1,2,\ldots,n$. 

\item At time index $t$, perform a measurement of $V_{bus}[t]$ and $I_{bus}[t]$ and compute the external impedance 
\begin{equation}
Z_l[t] = \frac{V_{bus}[t]}{I_{bus}[t]}
\label{eq:Zlm}
\end{equation}
and update the impedance matrix $D[t]$ at time index $t$ using (\ref{eq:impedancematrix}). 

\item Before the subsequent time step $t+1$, compute the module voltages $V_{opt}[t+1]$ according to 
\begin{equation}
V_{opt}[t+1] = \beta_{opt}[t] D[t] \left [ \begin{array}{cccc} \beta_1 & \beta_2 & \cdots & \beta_n \end{array} \right ]^T
\label{eq:Voptmp1}
\end{equation}
where $\beta_{opt}[t]$ is found by the LP problem in (\ref{eq:linprog}) using the updated impedance matrix $D[t]$ and communicate the $n$ elements $V_k[t+1]$ of $V_{opt}[t+1] = [V_1[t+1]~V_2[t+1]~ \cdots~ V_n[t+1]]^T$ to each of the modules $k=1,2,\ldots,n$.

\item At time step $t+1$, {\em all\/} modules $k=1,2,\ldots,n$ update the module voltage $V_k$ to $V_k = V_k[t+1]$.

\item Increment time index $t=t+1$ and restart at step $1)$.

\end{enumerate}

It should be noted that the recursive updates of $V_{opt}[t]$ explained above converges in a single time step in case $Z_l$ is fixed. In order to be able to track (fast) time-varying changes in the external load $Z_l$, the above procedure should allow high frequent measurements (and communication) of bus voltage $V_{bus}$ {\em and\/} bus current $I_{bus}$. Furthermore, the LP problem in (\ref{eq:linprog}) is solved and the results are used to communicate and update all the module voltages $V_{opt}[t+1] = [V_1~V_2~ \cdots V_n]^T$ in (\ref{eq:Voptmp1}). Since the LP program is solved and the elements of $V_{opt}[t+1]$ are communicated to each module, this is a centralized implementation of the recursive module scheduling. Each module $k=1,2,\ldots,n$ simply only receives its $V_k[t+1]$ from the (centrally) computed optimal LP solution $V_{opt}[t+1]$.

\section{EXPERIMENTAL SETUP}
\begin{figure}[ht]
\begin{center}
\includegraphics[width=0.9\columnwidth]{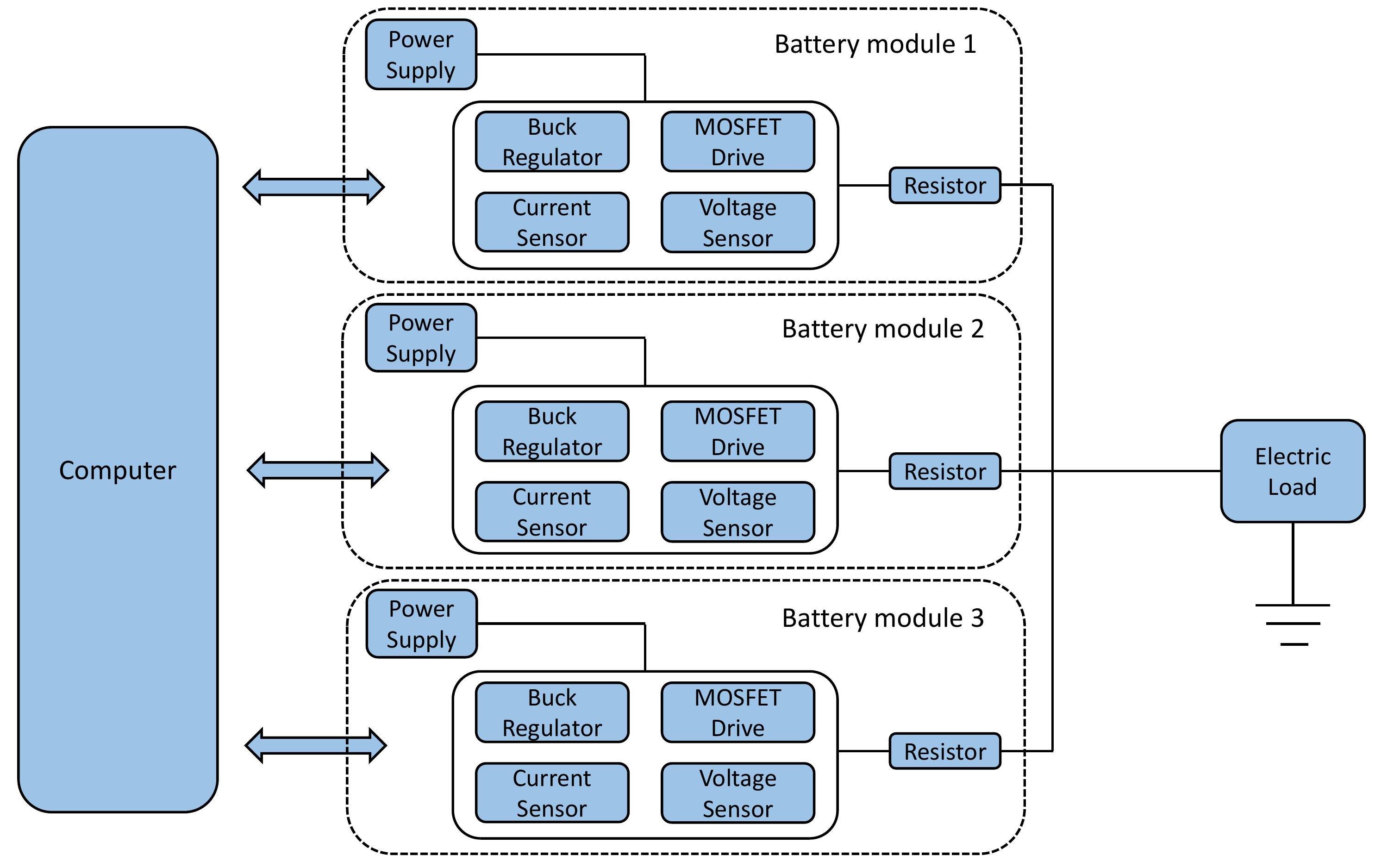}   
\caption{Schematic of the experimental battery tester.} 
\label{fig:schematic of the experimental battery tester}
\end{center}
\end{figure}

An experimental setup is built to demonstrate current scheduling, where the modulation demand signal can be applied and at the same time, voltage and current of each battery module can be measured. To explain the experimental setup, please refer to schematic diagram of Fig.~\ref{fig:schematic of the experimental battery tester} that indicates the parallel connection of 3 buck regulated battery modules. In addition, the parallel connection of 3 buck regulated battery modules is connected to a common DC bus and connected to an electrical load which is composed by a parallel connection of load resistors. Specifically, each parallel buck regulated battery module is composed of an adjustable power supply in series with an impedance, and a buck regulator. The buck regulator is composed by a pulse-width modulation (PWM) driven metal-oxide-semiconductor field-effect transistor (MOSFET), a fly-by diode, an inductor, and an Arduino Uno board. 

\begin{figure}[ht]
\begin{center}
\includegraphics[width=0.9\columnwidth]{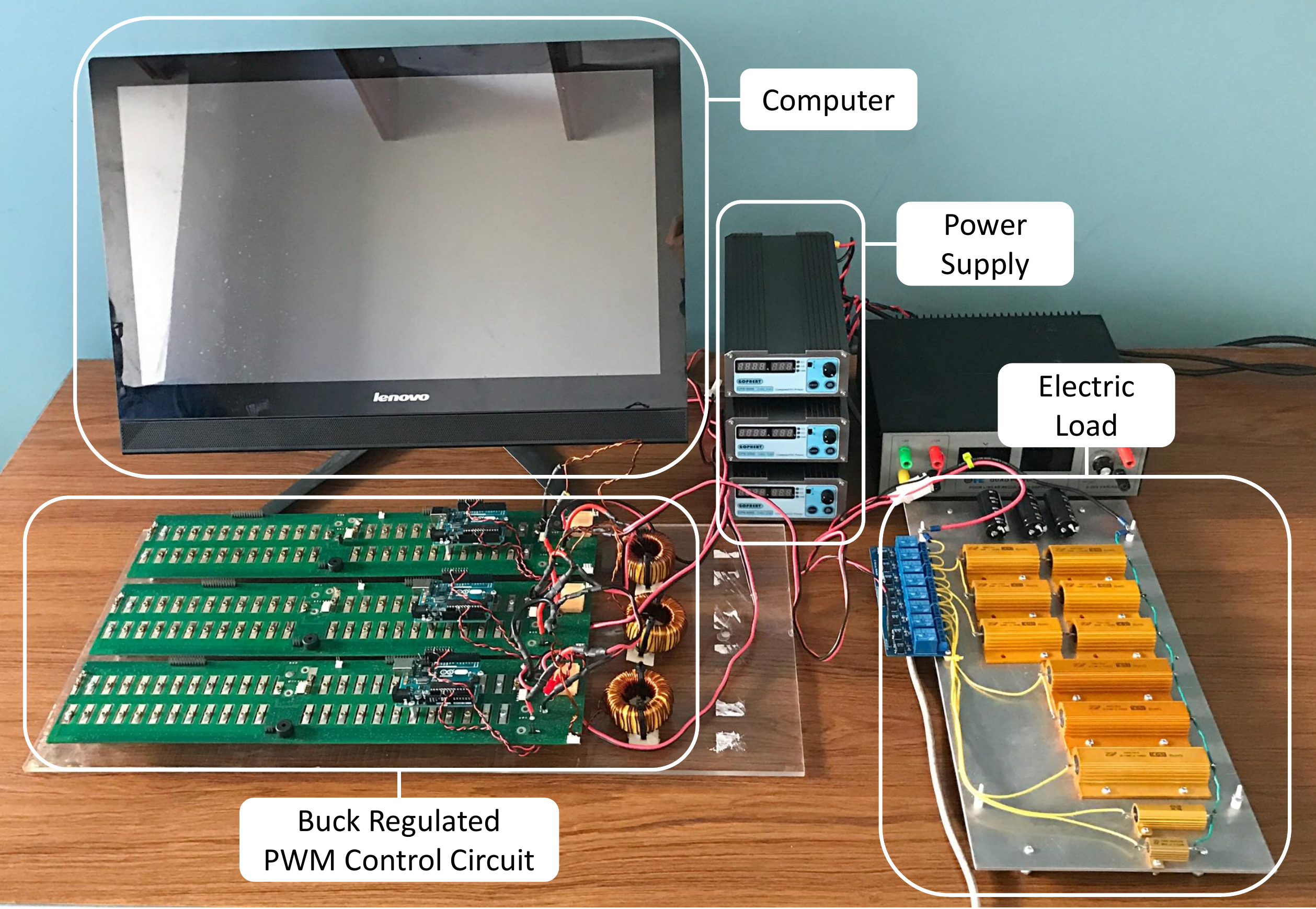}   
\caption{Photograph of the experimental battery tester.} 
\label{fig:battery tester}
\end{center}
\end{figure}

The experimental setup tester description is summarized by a photograph shown in Fig.~\ref{fig:battery tester}. The power supply used in the experiment is GOPHERT CPS-6005 0-60V 0-5A Adjustable Switched Mode DC Power Supply. The MOSFET on the buck regulator is switched by corresponding control signals sent from PWM pins of the Arduino Uno board to modulate down OCV of that battery module. The Arduino Uno board can be employed to get current/voltage measured real-time signals by its analog input pins, and can also communicate with the computer through USB cable. In the computer, the MATLAB-Arduino interface is applied to automatically implement current scheduling algorithms and save measured real-time data simultaneously. Moreover, the MOSFETs applied in the experiment are driven by 62.5 kHz PWM modulation output frequency from Arduino PWM pins, and they are with low drain-to-source on-resistance that is suitable for high current of battery modules.

\section{EXPERIMENTAL RESULTS }

\begin{figure}
\begin{center}
\includegraphics[width=\columnwidth]{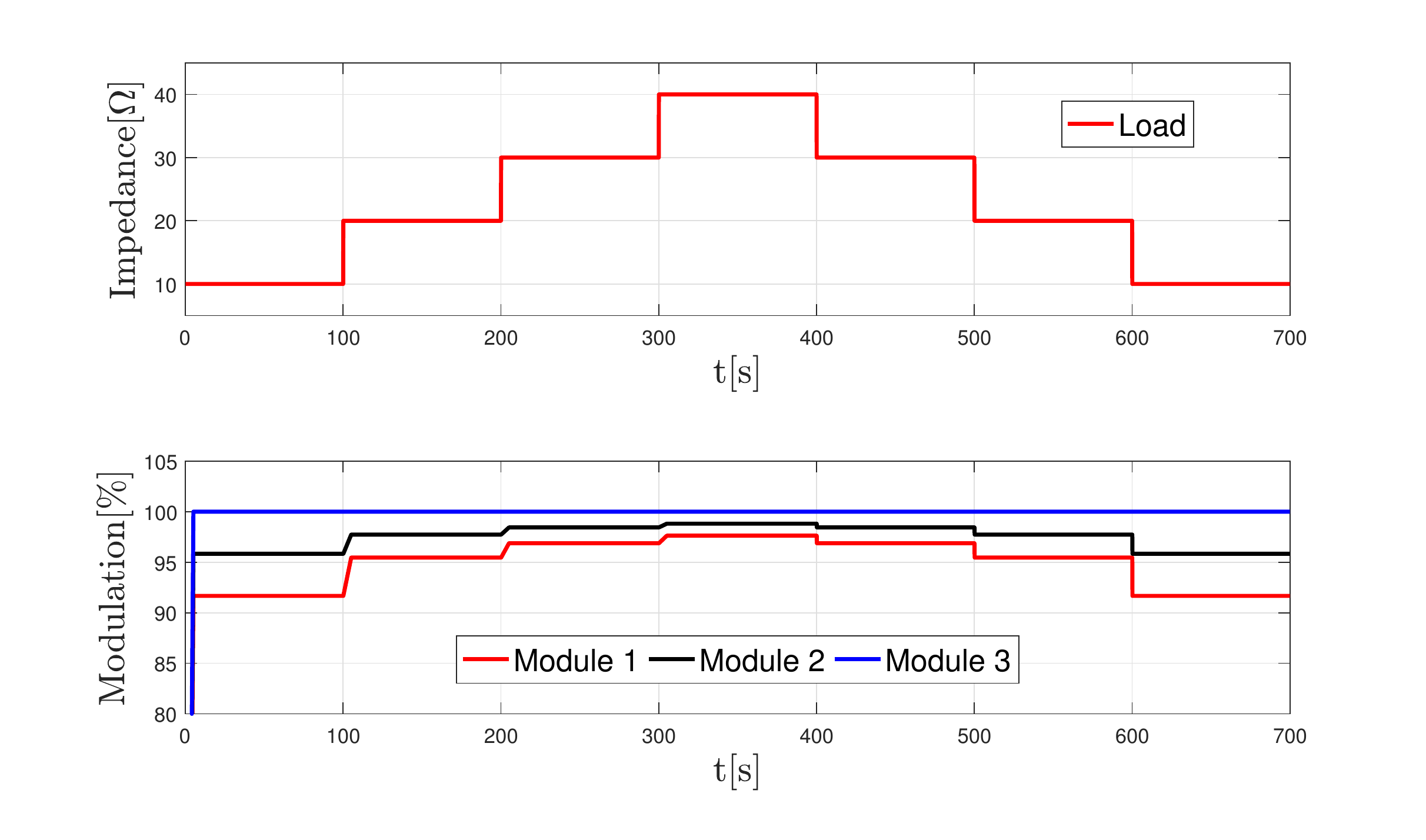}   
\caption{Time-varying external load impedance (top figure), and modulation of battery voltage for recursive centralized module scheduling of 3 parallel placed battery modules (bottom figure).} 
\label{fig:load impedance and modulation experiment results}
\end{center}
\end{figure}

For experimental verification of the centralized recursive balanced scheduling, 3 parallel placed battery modules with scaled down OCVs $V_1^{OCV}$ = $V_2^{OCV}$ = $V_3^{OCV}$ = 5V and internal impedance values $Z_1$ = 3$\Omega$, $Z_2$ = 4.5$\Omega$, and $Z_3$ = 6$\Omega$ are used, subjected to a time-varying external load shown at the top of the Fig.~\ref{fig:load impedance and modulation experiment results}. Specifically, in entire 700s experimental period, the external load is automatically increased 10 $\Omega$ every 100s from $t_0$ = 0s to $t_1$ = 400s, and decreased 10 $\Omega$ every 100s from $t_1$ = 400s to $t_2$ = 700s. 

The resulting time-varying PWM modulation factors for the recursive updates of the internal voltage is shown in the bottom of the Fig.~\ref{fig:load impedance and modulation experiment results}. Module 3 is always fixed at 100$\%$ modulation, because of its highest impedance $Z_3$ = 6$\Omega$, which needs to modulate down other battery modules to keep equal (balanced) currents. It should be noted that increasing PWM requires ramp function with 5s ramp-up period and decreasing PWM can happen instantaneously in order to protect the battery modules.

\begin{figure}
\begin{center}
\includegraphics[width=\columnwidth]{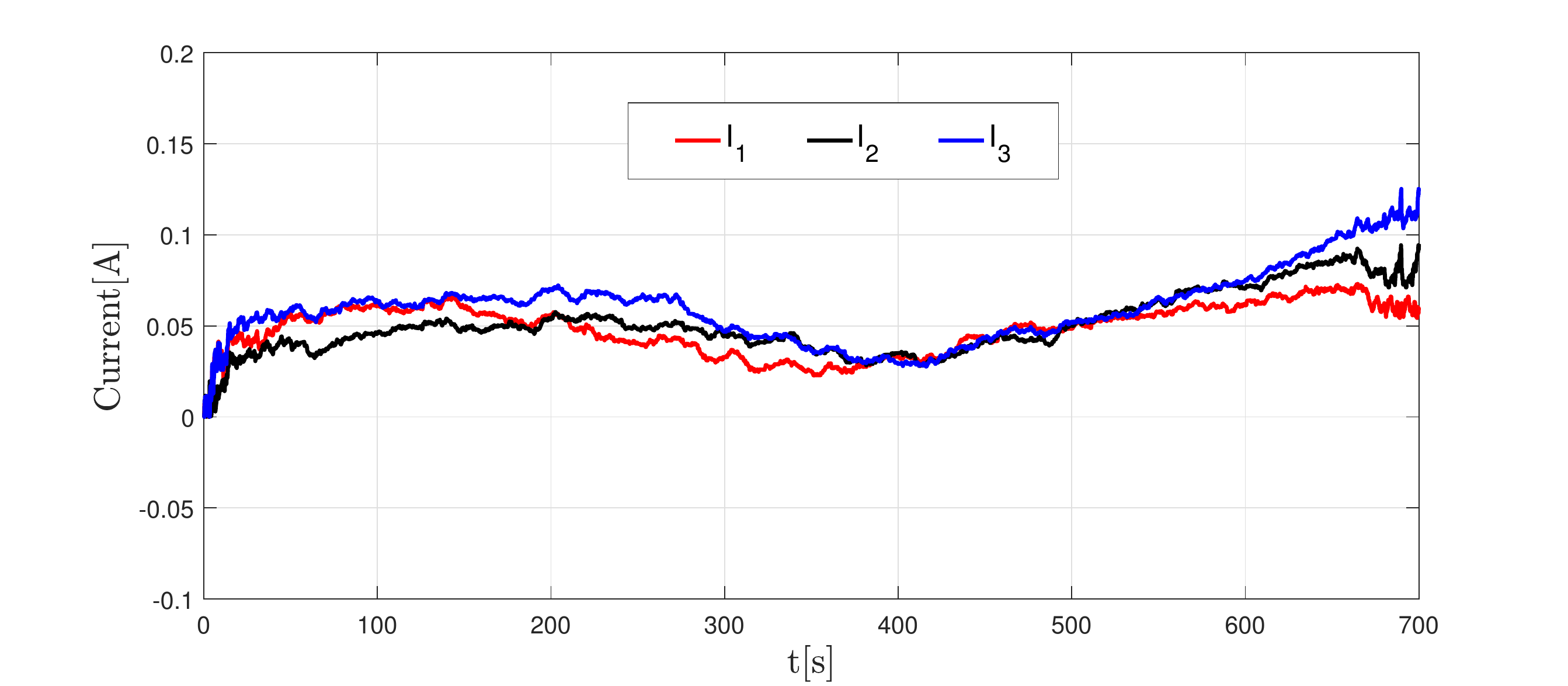}   
\caption{Current in battery modules for recursive centralized module scheduling of 3 parallel placed battery modules.} 
\label{fig:current measurements}
\end{center}
\end{figure}

The experimental results for centralized recursive time-varying balanced scheduling are shown in Fig.~\ref{fig:current measurements}. It can be seen that individual module currents keep relatively closed to each other in spite of time-varying external load, which verifies the feasibility of proposed centralized recursive scheduling for balancing individual battery module current. 

\section{CONCLUSIONS}

The efficiency and flexibility of a battery system that consists of parallel placed battery modules can be significantly improved when partly empty or failing modules of the battery pack can be exchanged for fast charging and fault correction capabilities. However, control or scheduling of battery modules is required to account for differences between state of charge, instantaneous and nominal capacity, and internal impedances of the battery modules.

In this paper a solution is provided that allows for centralized recursive current scheduling of parallel placed battery modules. The current scheduling algorithm uses a Linear Programming formulation to compute the optimal open circuit voltage values of each battery module so that currents of all battery modules are balanced to avoid stray currents between modules. This centralized current scheduling can adjust individual module currents to be equal (balanced) and increase the total bus current after optimization. Furthermore, an experimental setup with 3 parallel battery modules validates the optimal current scheduling algorithms. The experimental results indicate that the proposed method is able to effectively balance (equal) individual battery currents under time-varying external load conditions. 

The future work of this study is to propose a decentralized recursive current scheduling method by solving the similar LP problem to efficiently reduce the centralized communication requirements on speed and reliability of the communication hardwares. By doing so, we hope to show that the decentralized solution can balance (equal) individual module currents and optimize total bus current in order to eliminate the need for high speed central communication between battery modules.

\addtolength{\textheight}{-12cm}   





\bibliographystyle{IEEEtran}
\bibliography{Reference}

\end{document}